\documentclass[11pt]{article}
\usepackage{epic,latexsym,amssymb}
\usepackage{amsfonts}
\usepackage{amscd}
\usepackage{amsmath}
\usepackage{graphicx}
\usepackage{color}
\usepackage{caption,
caption}
\usepackage{tikz}

\usepackage{float}

\newtheorem{thm}{Theorem}

\newtheorem{ob}[thm]{Observation}
\newtheorem{prop}[thm]{Proposition}
\newtheorem{lem}[thm]{Lemma}

\newtheorem{claim}{Claim}

\newcommand{\trim}{\mathrm{trim}}

\newcommand{\diam}{{\rm diam}}

\newcommand{\cF}{\mathcal{F}}

\newcommand{\cT}{\mathcal{T}}
\newcommand{\cO}{\mathcal{O}}
\newcommand{\cH}{\mathcal{H}}

\newcommand{\proof}{\noindent\textbf{Proof. }}

\newcommand{\smallqed}{{\tiny ($\Box$)}}
\newcommand{\qed}{$\Box$}

\newcommand{\1}{\vspace{0.1cm}}
\newcommand{\2}{\vspace{0.2cm}}

\newcommand{\QEDmark}{\mbox{\textsc{qed}}}
\newcommand{\proofStarter}[1]{\textsc{#1} }

\def\vertex(#1){\put(#1){\circle*{2}}}
\def\vertexo(#1){\put(#1){\circle{2}}}
\def\vert(#1){\put(#1){\circle*{1.5}}}
\def\verto(#1){\put(#1){\circle{1.5}}}
\def\lab(#1)#2{\put(#1){\makebox(0,0)[c]{#2}}}
\definecolor{DarkGreen}{rgb}{0.2, 0.6, 0.3}

\definecolor{electricindigo}{rgb}{0.44, 0.0, 1.0}

\let\oldenumerate\enumerate
\renewcommand{\enumerate}{
  \oldenumerate
  \setlength{\itemsep}{0.5pt}
  \setlength{\parskip}{0pt}
  \setlength{\parsep}{0pt}
}

\begin{document}

\title{Total Forcing Sets in Trees}
\author{$^{1,2}$Randy Davila and $^{1}$Michael A. Henning\\
\\
$^1$Department of Pure and Applied Mathematics\\
University of Johannesburg \\
Auckland Park 2006, South Africa \\
\\
$^2$Department of Mathematics \\
Texas State University \\
San Marcos, TX 78666, USA \\
\small {\tt Email: rrd32@txstate.edu}
}

\date{}
\maketitle

\begin{abstract}
A dynamic coloring of the vertices of a graph $G$ starts with an initial subset $S$ of colored vertices, with all remaining vertices being non-colored. At each discrete time interval, a colored vertex with exactly one non-colored neighbor forces this non-colored neighbor to be colored. The initial set $S$ is called a forcing set of $G$ if, by iteratively applying the forcing process, every vertex in $G$ becomes colored. If the initial set $S$ has the added property that it induces a subgraph of $G$ without isolated vertices, then $S$ is called a total forcing set in $G$. The minimum cardinality of a total forcing set in $G$ is its total forcing number, denoted $F_t(G)$. We prove that if $T$ is a tree of order $n \ge 3$ with maximum degree~$\Delta$, then $F_t(T) \le \frac{1}{\Delta}((\Delta - 1)n + 1)$, and we characterize the infinite family of trees achieving equality in this bound. We also prove that if $T$ is a non-trivial tree with $n_1$ leaves, then $F_t(T) \ge n_1$, and we characterize the infinite family of trees achieving equality in this bound. As a consequence of this result, the total forcing number of a non-trivial tree is strictly greater than its forcing number. In particular, we prove that if $T$ is a non-trivial tree, then $F_t(T) \ge F(T)+1$, and we characterize extremal trees achieving this bound.
\end{abstract}


\newpage
\section{Introduction}

A dynamic coloring of the vertices in a graph is a coloring of the vertex set which may change, or propagate, throughout the vertices during discrete time intervals. Of the dynamic colorings, the notion of \emph{forcing sets} (\emph{zero forcing sets}), and the associated graph invariant known as the \emph{forcing number} (\emph{zero forcing number}), are arguably the most prominent, see for example \cite{AIM-Workshop, zf_np, Davila Kenter, Genter1, Genter2, Edholm, zf_np2}. 
In the study of minimum forcing sets in graphs, it is natural to consider the initial structure of such sets. In particular, if a forcing set induces an isolate-free subgraph, then the set in question is called a \emph{total forcing set}.

More formally, let $G$ be a graph with vertex set $V = V(G)$ and edge set $E = E(G)$. The \emph{forcing process} is defined in~\cite{DaHe17+} as follows: Let $S \subseteq V$ be a set of initially ``colored" vertices, all other vertices are said to be ``non-colored". A vertex contained in $S$ is said to be $S$-colored, while a vertex not in $S$ is said to be $S$-uncolored. At each time step, if a colored vertex has exactly one non-colored neighbor, then this colored vertex \emph{forces} its non-colored neighbor to become colored. If $v$ is such a colored vertex, we say that $v$ is a \emph{forcing vertex}. We say that $S$ is a \emph{forcing set}, if by iteratively applying the forcing process, all of $V$ becomes colored. We call such a set $S$, an $S$-forcing set. In addition, if $S$ is an $S$-forcing set in G and $v$ is a $S$-colored vertex that forces a new vertex to be colored, then $v$ is an \emph{$S$-forcing vertex}. The cardinality of a minimum forcing set in $G$ is the \emph{forcing number} of $G$, denoted $F(G)$.

If $S$ is a forcing set which also induces a graph without isolated vertices, then $S$ is a \emph{total forcing set}, abbreviated as a TF-\emph{set} of $G$. The \emph{total forcing number} of $G$, written $F_t(G)$, is the cardinality of a minimum TF-set in $G$. The concept of a total forcing set was first introduced and studied by Davila in~\cite{Davila}, and is further studied, for example, by the authors in~\cite{DaHe17+}. In this paper, we study total forcing sets in trees. In particular, we study trees with smallest possible total forcing number, as well as trees with largest possible forcing number.

\medskip
\noindent\textbf{Definitions and Notation.}
For notation and graph terminology, we will typically follow \cite{MHAYbookTD}. Throughout this paper, all graphs will be considered undirected, simple and finite. Specifically, let $G$ be a graph with vertex set $V(G)$ and edge set $E(G)$, and of order $n = |V(G|$ and size $m = |E(G)|$. If the graph $G$ is clear from the context, we simply write $V$ and $E$ rather than $V(G)$ and $E(G)$, and we write $G = (V,E)$. A \emph{non}-\emph{trivial} graph is a graph with at least two vertices.

Two vertices $v$ and $w$ are adjacent, or neighbors,  in $G$ if $vw \in E$. The open neighborhood of a vertex $v \in V$, is the set of neighbors of $v$, denoted $N_G(v)$, whereas its closed neighborhood is $N_G[v] = N_G(v) \cup \{v\}$. The open neighborhood of $S \subseteq V$ is the set of all neighbors of vertices in $S$, denoted $N_G(S)$, whereas the closed neighborhood of $S$ is $N_G[S] = N_G(S) \cup S$. We denote the \emph{degree} of a vertex $v$ in a graph $G$ by $d_G(v)$, or simply by $d(v)$ if the graph $G$ is clear from the context. Thus, $d_G(v) = |N_G(v)|$. The minimum and maximum degree among the vertices of $G$ is denoted by $\delta = \delta(G)$ and $\Delta = \Delta(G)$, respectively.
For a set of vertices $S \subseteq V$, the subgraph induced by $S$ is denoted by $G[S]$. The subgraph obtained from $G$ by deleting all vertices in $S$ and all edges incident with vertices in $S$ is denoted by $G - S$. If $S = \{v\}$, we simply write $G - v$ rather than $G - S$.

A \emph{leaf} is a vertex of degree~$1$, while its neighbor is a \emph{support vertex}. A \emph{strong support vertex} is a vertex with at least two leaf neighbors. A \emph{star} is a non-trivial tree with at most one vertex that is not a leaf. Thus, a star is the tree $K_{1,k}$ for some $k \ge 1$. For $r, s \ge 1$, a \emph{double star} $S(r,s)$ is the tree with exactly two vertices that are not leaves, one of which has $r$ leaf neighbors and the other $s$ leaf neighbors. We will denote a path on $n$ vertices by $P_n$. We define a \emph{pendant edge} of a graph to be an edge incident with a vertex of degree~$1$.

The distance between two vertices $v$ and $w$ in  $G$ is the length of a shortest $(v,w)$-path in $G$, and is denoted by $d_G(v,w)$. If no $(v,w)$-path exists in $G$, then we define $d_G(v,w) =\infty$. The maximum distance among all pairs of vertices of $G$ is the \emph{diameter} of $G$, denoted by $\diam(G)$. The \emph{eccentricity} of a vertex $v$ in $G$ is the maximum distance of a vertex from $v$ in~$G$. A vertex of minimum eccentricity is called a \emph{central vertex} of $G$. In particular, the central vertex of a star of order at least~$3$ is the vertex that is not a leaf, while a double star contains two central vertices, namely the two vertices that are not leaves.

A \emph{rooted tree} $T$ distinguishes one vertex $r$ called the \emph{root}. For each vertex $v \ne r$ of $T$, the \emph{parent} of $v$ is the
neighbor of $v$ on the unique $(r,v)$-path, while a \emph{child} of $v$ is any other neighbor of $v$. The set of children of $v$ is denoted by  $C(v)$. A \emph{descendant} of $v$ is a vertex $u \ne v$ such that the unique $(r,u)$-path contains $v$, while an \emph{ancestor} of $v$ is a vertex $u \ne v$ that belongs to the $(r,v)$-path in $T$. In particular, every child of $v$ is a descendant of~$v$ while the parent of $v$ is an ancestor of $v$. The \emph{grandparent} of $v$ is the ancestor of $v$ at distance~$2$ from $v$. A \emph{grandchild} of $v$ is the descendant of $v$ at distance~$2$ from $v$.
We let $D(v)$ denote the set of descendants of $v$, and we define $D[v] = D(v) \cup \{v\}$. The \emph{maximal subtree} at $v$ is the subtree of $T$ induced by $D[v]$, and is denoted by $T_v$.

We use the standard notation $[k] = \{1,\ldots,k\}$.

\section{Main Results}

Recently, the authors established the following upper bound on the total forcing number of a graph with minimum degree at least two in terms of the order and maximum degree of the graph.

\begin{thm}{\rm (\cite{DaHe17+})}
\label{t:upperbd}
If $G$ is a connected graph of order $n \ge 3$  with minimum degree at least~$2$ and with maximum degree~$\Delta$, then
\[
F_t(G) \le \left( \frac{\Delta}{\Delta + 1} \right) n,
\]
with equality if and only if $G \cong K_n$. \end{thm}

In this paper, we study the total forcing number of a tree. We have four immediate aims. First to prove that the upper bound established in Theorem~\ref{t:upperbd} also holds for the class of trees. Secondly, to establish a much stronger upper bound on the total forcing number of a tree in terms of its order and maximum degree, and to characterize the extremal trees. Thirdly, to establish a lower bound on the total forcing number of a tree in terms of the number of leaves in the tree, and once again to characterize the extremal trees. More precisely, we shall prove the following four results, where $\cT$, $\cF$, and $\cH$ are families of trees we construct in Section~\ref{S:cT},~\ref{S:cF}, and~\ref{S:cS} respectively.

\begin{thm}
\label{t:tree0}
If $T$ is a tree of order~$n \ge 3$ with maximum degree~$\Delta$, then
\[
F_t(T) \le \left( \frac{\Delta}{\Delta + 1} \right) n,
\]
with equality if and only if $T \cong K_{1,\Delta}$.
\end{thm}

\begin{thm}
\label{t:tree1}
If $T$ is a tree of order~$n \ge 3$ with maximum degree~$\Delta$, then
\[
F_t(T) \le \frac{(\Delta - 1)n + 1}{\Delta},
\]
with equality if and only if $T \in \cT_{\Delta}$.
\end{thm}

\begin{thm}\label{t:tree2}
If $T$ is a non-trivial tree with $n_1$ leaves, then
$F_t(T) \ge n_1$,
with equality if and only if $T\in \cF$.
\end{thm}

\begin{thm}\label{t:tree3}
If $T$ is a non-trivial tree, then
$F_t(T)\ge F(T) + 1$,
with equality if and only if $T\in \cH$.
\end{thm}

\section{Known Results}

Before proceeding with a proof of our main results, we present some known results first observed in~\cite{DaHe17+}.

\begin{ob}{\rm (\cite{DaHe17+})}
\label{ob:1}
If $G$ is an isolate-free graph of order $n\ge 3$, then $F_t(G)= n - 1$, with equality if and only if $G = K_n$ or $G = K_{1,n-1}$.
\end{ob}

\begin{ob}{\rm (\cite{DaHe17+})}
\label{ob:2}
Every total forcing set in an isolate-free graph contains every strong support vertex of the graph and all except possibly one leaf neighbor of each strong support vertex.
\end{ob}


\section{Proof of Theorem~\ref{t:tree1}}
\label{S:cT}

In this section, we present a proof of Theorem~\ref{t:tree1}. In order to determine an upper bound on the total forcing number of a tree, we define a family $\cT_{\Delta}$ of trees as follows.

\medskip
\noindent\textbf{The family $\cT_{\Delta}$.} Let $\cT_{\Delta}$ be the family of all trees $T$ with maximum degree~$\Delta$ whose vertex set $V(T)$ can be partitioned into sets $(V_1,\ldots,V_k)$ such that the following holds, where $T_i = T[V_i]$ for $i \in [k]$.
\\[-20pt]
\begin{enumerate}
\item[$\bullet$] $T_1 \cong K_{1,\Delta}$, and if $k \ge 2$, then $T_i \cong K_{1,\Delta - 1}$  for $i \in [k] \setminus \{1\}$.
\item[$\bullet$] For $i \in [k]$, the central vertex $v_i$ of the star $T_i$ is a strong support vertex of degree~$\Delta$ in the tree $T$.
\item[$\bullet$] The set $\{v_1,\ldots,v_k\}$ is an independent set in $T$.
\end{enumerate}

\begin{figure}[htb]
\begin{center}

\begin{tikzpicture}[scale=.8,style=thick,x=1cm,y=1cm]
\def\vr{2.5pt} 
\path (0,0) coordinate (u1);
\path (0.75,0) coordinate (u2);
\path (1.5,0) coordinate (u3);
\path (0.75,1) coordinate (v);
%
\draw (v) -- (u1);
\draw (v) -- (u2);
\draw (v) -- (u3);
\draw (v) [fill=black] circle (\vr);
\draw (u1) [fill=black] circle (\vr);
\draw (u2) [fill=black] circle (\vr);
\draw (u3) [fill=black] circle (\vr);
\draw (0.75,-0.75) node {(a)};
\draw[anchor = south] (v) node {{\small $v_1$}};
\path (5,0) coordinate (u1);
\path (5.5,1) coordinate (v1);
\path (6,0) coordinate (u2);
\path (6.5,2) coordinate (v);
\path (7,0) coordinate (u3);
\path (7.5,1) coordinate (v2);
\path (8,0) coordinate (u4);
\draw (v) -- (v1);
\draw (v) -- (v2);
\draw (v1) -- (u1);
\draw (v1) -- (u2);
\draw (v2) -- (u3);
\draw (v2) -- (u4);
\draw (v) [fill=black] circle (\vr);
\draw (v1) [fill=black] circle (\vr);
\draw (v2) [fill=black] circle (\vr);
\draw (u1) [fill=black] circle (\vr);
\draw (u2) [fill=black] circle (\vr);
\draw (u3) [fill=black] circle (\vr);
\draw (u4) [fill=black] circle (\vr);
\draw (6.5,-0.75) node {(b)};
\draw[anchor = east] (v1) node {{\small $v_1$}};
\draw[anchor = west] (v2) node {{\small $v_2$}};
\path (11.5,0) coordinate (u1);
\path (12.5,0) coordinate (u2);
\path (13.5,0) coordinate (u3);
\path (14.5,0) coordinate (u4);
\path (15.5,0) coordinate (u5);
\path (16.5,0) coordinate (u6);
\path (14,2) coordinate (v);
\path (12,1) coordinate (v1);
\path (14,1) coordinate (v2);
\path (16,1) coordinate (v3);
\draw (v) -- (v1);
\draw (v) -- (v2);
\draw (v) -- (v3);
\draw (v1) -- (u1);
\draw (v1) -- (u2);
\draw (v2) -- (u3);
\draw (v2) -- (u4);
\draw (v3) -- (u5);
\draw (v3) -- (u6);
\draw (v) [fill=black] circle (\vr);
\draw (v1) [fill=black] circle (\vr);
\draw (v2) [fill=black] circle (\vr);
\draw (v3) [fill=black] circle (\vr);
\draw (u1) [fill=black] circle (\vr);
\draw (u2) [fill=black] circle (\vr);
\draw (u3) [fill=black] circle (\vr);
\draw (u4) [fill=black] circle (\vr);
\draw (u5) [fill=black] circle (\vr);
\draw (u6) [fill=black] circle (\vr);
\draw (14,-0.75) node {(c)};
\draw[anchor = east] (v1) node {{\small $v_1$}};
\draw[anchor = east] (v2) node {{\small $v_2$}};
\draw[anchor = west] (v3) node {{\small $v_3$}};
\end{tikzpicture}
\end{center}
\vskip -0.4 cm
\caption{The three trees in the family $\cT_3$.} \label{f:cT}
\end{figure}
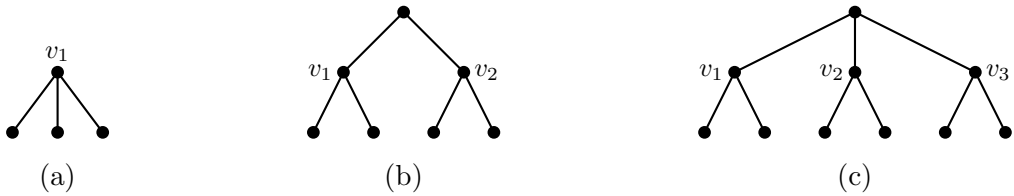

We call the trees $T_1,\ldots,T_k$ the \emph{underlying subtrees} of the tree $T$.
We note that $\cT_2$ consists only of the path $P_3$; that is, $\cT_2 = \{P_3\}$. The family $\cT_3$ consists of the three trees shown in Figure~\ref{f:cT}(a),~\ref{f:cT}(b) and~\ref{f:cT}(c).

\bigskip

We first establish useful properties of trees in the family~$\cT_{\Delta}$.

\begin{lem}
\label{l:lem1}
If $T$ is a tree of order~$n$
that belongs to the family $\cT_\Delta$, then the following holds. \\[-27pt]
\begin{enumerate}
\item $F_t(T) = \frac{1}{\Delta}((\Delta - 1)n + 1)$.
\item The set consisting of all vertices of $T$, except for exactly one leaf neighbor in $T$ of the central vertex of each underlying subtree of $T$ is a minimum TF-set of $T$.
\end{enumerate}
\end{lem}
\proof Let $T \in \cT_\Delta$ be a tree of order~$n$ with maximum degree~$\Delta$. We proceed by induction on the number, $k$, of underlying subtrees of the tree $T$. If $k = 1$, then $T \cong K_{1,\Delta}$ and by Observation~\ref{ob:1}, $F_t(T) = \Delta = ((\Delta - 1)n + 1)/\Delta$ noting that here $n = \Delta + 1$. Further, Property~(b) is immediate in this case. This establishes the base case. Let $k \ge 2$ and assume that if $T'$ is a tree of order~$n'$ in the family $\cT_\Delta$ with $k'$ underlying subtrees, where $k' < k$, then $F_t(T') = \frac{1}{\Delta}((\Delta - 1)n' + 1)$.

Let $T$ be a tree of order~$n$ in the family $\cT_\Delta$ with $k$ underlying subtrees given by $T_1,T_2,\ldots,T_k$, where $T_1 \cong K_{1,\Delta}$ and $T_i \cong K_{1,\Delta - 1}$ for $i \in [k] \setminus \{1\}$. Recall that $V_i = V(T_i)$ for $i \in [k]$. We note that $n = k \Delta + 1$. Let $S$ consist of all vertices of $T$, except for exactly one leaf neighbor of the central vertex of each subtree of $T$.
Further, let $S_i = S \cap V_i$ for $i \in [k]$. Thus, $|S_1| = \Delta$ and $|S_i| = \Delta - 1$ for $i \in [k] \setminus \{1\}$. The set $S$ is a TF-set, and so
\1
\[
\begin{array}{lcl}
F_t(T) \, \le  \, |S| & = & \displaystyle{ \sum_{i=1}^k |S_i| } \1 \\
& = & \displaystyle{ \Delta  + \sum_{i=2}^k (\Delta - 1) } \2 \\
& = & k(\Delta - 1) + 1 \2 \\
& = & \frac{1}{\Delta}((\Delta - 1)n + 1).
\end{array}
\]

If $F_t(T) = \frac{1}{\Delta}((\Delta - 1)n + 1)$, then we must have equality throughout the above inequality chain, implying that $F_t(T) = |S|$ and hence that $S$ is a minimum TF-set of~$T$. Therefore, if Property~(a) holds, then Property~(b) holds. Hence, it suffices for us to prove Property~(a); that is, to prove that $F_t(T) = \frac{1}{\Delta}((\Delta - 1)n + 1)$.

Let $T^*$ be the graph of order~$k$ whose vertices correspond to the $k$ subtrees of $T$, and where we add an edge between two vertices of $T^*$ if the corresponding subtrees of $T$ are joined by an edge in $T$. Since $T$ is a tree, so too is $T^*$. Since $k \ge 2$, at least one leaf of $T^*$ corresponds to an underlying subtree $T_i$ for some $i \ge 2$. Renaming the subtrees $T_2,\ldots,T_k$, if necessary, we may assume that the underlying subtree $T_k$ corresponds with such a leaf of $T^*$. Since $v_k$ has degree $\Delta$ in $T$ and degree~$\Delta - 1$ in $T_k$, this implies that the $\Delta - 1$ leaf neighbors of $v_k$ in $T_k$ are all leaf neighbors of $v_k$ in $T$ and that the vertex $v_k$ is joined to exactly one vertex, say $x$, in $T$ that does not belong to $V(T_k)$. Let $x$ belong to the underlying subtree $T_j$, where $j \in [k-1]$. By definition of the family $\cT_\Delta$, the set of central vertices of the $k$ underlying subtrees form an independent set in $T$, implying that $x \ne v_j$.

Let $T'$ be the tree of order~$n'$ obtained from $T$ by deleting the vertices in $V(T_k)$. We note that $T' \in \cT_\Delta$ with underlying subtrees $T_1,\ldots,T_{k-1}$. Further, $n' = n - \Delta$. Applying the inductive hypothesis to the tree $T'$, $F_t(T') = \frac{1}{\Delta}((\Delta - 1)n' + 1)$. Let $v_k'$ be an arbitrary leaf neighbor of $v_k$ in $T_k$ (and therefore in $T$). If $S'$ is a minimum TF-set of $T'$, then $S' \cup (V(T_k) \setminus \{v_k'\})$ is a TF-set of $T'$, implying that $F_t(T) \le |S'| + \Delta - 1 = F_t(T') + \Delta - 1$. Conversely, let $X$ be a minimum TF-set in $T$. By Observation~\ref{ob:2}, the set $X$ contains the vertex $v_k$ and all except possibly one leaf neighbor of $v_k$. If the set $X$ contains all leaf neighbors of $v_k$, then by the minimality of the set $X$, we note that $x \notin X$. However, in this case, we can simply remove exactly one leaf neighbor of $v_k$ from $X$ and add the vertex $x$ to the set $X$ to produce a new minimum TF-set in $T$. Thus, we may choose the set $X$ so that $x \in X$ and $v_k' \notin X$. Let $X'$ be the restriction of $X$ to $V(T')$; that is, $X' = X \cap V(T')$. Since $X$ is a TF-set of $T$, the set $X'$ is a TF-set of $T'$. Thus, $F_t(T') \le |X'| = |X| - (\Delta - 1) = F_t(T) - \Delta + 1$. Consequently, \1
\[
\begin{array}{lcl}
F_t(T) & = & F_t(T') + \Delta - 1 \1 \\
& = & \frac{1}{\Delta}((\Delta - 1)n' + 1) + \Delta - 1 \2 \\
& = & \frac{1}{\Delta}((\Delta - 1)(n - \Delta) + 1) + \Delta - 1 \2 \\
& = & \frac{1}{\Delta}((\Delta - 1)n + 1). 
\end{array}
\]
This completes the proof of Lemma~\ref{l:lem1}.~\qed

\medskip
We are now in a position to prove Theorem~\ref{t:tree1}. Recall its statement.

\noindent \textbf{Theorem~\ref{t:tree1}}. \emph{If $T$ is a tree of order~$n \ge 3$ with maximum degree~$\Delta$, then
\[
F_t(T) \le \frac{(\Delta - 1)n + 1}{\Delta},
\]
with equality if and only if $T \in \cT_{\Delta}$.
}

\medskip
\proof  We proceed by induction on the order~$n \ge 3$ of a tree $T$ with maximum degree~$\Delta$. If $n = 3$, then $T \cong P_3 \in \cT_2$, $\Delta = 2$ and $F_t(T) = 2 = ((\Delta - 1)n + 1)/\Delta$. This establishes the base case. Let $n \ge 4$ and assume that if $T'$ is a tree of order~$n'$ and maximum degree~$\Delta'$, where $3 \le n' < n$ and $\Delta' \le \Delta$, then $F_t(T') \le ((\Delta' - 1)n' + 1)/\Delta'$, with equality if and only if $T' \in \cT_{\Delta'}$. If $\Delta' < \Delta$, then, by monotonicity, we note that $((\Delta' - 1)n' + 1)/\Delta' < ((\Delta - 1)n' + 1)/\Delta$. Let $T$ be a tree of order~$n$ and maximum degree~$\Delta$. We note that $\Delta \ge 2$.

Suppose that $\Delta = 2$, and so $T \cong P_n$ where we recall that $n \ge 4$. Since every non-trivial path has total forcing number 2, we observe that $F_t(T) = 2 < (n+1)/2 = ((\Delta - 1)n + 1)/\Delta$. Hence, we may assume that $\Delta \ge 3$, for otherwise the desired result follows.

Suppose that $\diam(T) = 2$, and so $T$ is a star. In this case, $\Delta = n-1$ and $T \cong K_{1,\Delta} \in \cT_{\Delta}$. By Observation~\ref{ob:1}, $F_t(T) = \Delta = ((\Delta - 1)n + 1)/\Delta$.  Hence, we may assume that $\diam(T) \ge 3$.

Suppose that $\diam(T) = 3$, and so $T \cong S(r,s)$ is a double star, where $1 \le r \le s$. Let $u$ and $v$ be the two vertices of $T$ that are not leaves, where $u$ has $r$ leaf neighbors and $v$ has $s$ leaf neighbors. Since $T$ has maximum degree~$\Delta$, we note that $s = \Delta - 1$, and so $n = r + s + 2 = r + \Delta + 1$. Let $u'$ and $v'$ be arbitrary leaf neighbors of $u$ and $v$, respectively. The set $V(T) \setminus \{u',v'\}$ is a TF-set of $T$, implying that $F_t(T) \le r + s = r + \Delta - 1$. Moreover, since $\Delta \ge r + 1$,
\1
\[
\begin{array}{lcl}
\displaystyle{ \frac{(\Delta - 1)n + 1}{\Delta} } & > &  \displaystyle{ \frac{(\Delta - 1)n}{\Delta}  } \2 \\
& = & \displaystyle{ \frac{1}{\Delta} ((\Delta - 1)(r + \Delta + 1)) } \1 \\
& = & \displaystyle{ r + \Delta - 1 + \frac{1}{\Delta} (\Delta - r - 1)} \1 \\
& \ge & r + \Delta - 1 \1 \\
& \ge & F_t(T). 
\end{array}
\]

Hence, we may assume that $\diam(T) \ge 4$, for otherwise the desired result holds. Let $u$ and $r$ be two vertices at maximum distance apart in $T$. Necessarily, $u$ and $r$ are leaves and $d(u,r) = \diam(T)$. We now root the tree $T$ at the vertex $r$. Let $v$ be the parent of $u$, $w$ the parent of $v$, $x$ be the parent of $w$, and $y$ the parent of~$x$.
We note that if $\diam(T) = 4$, then $y = r$; otherwise, $y \ne r$.

Let $d_T(v) = \ell$, where we note that $\ell \le \Delta$. Let $T'$ be the tree obtained from $T$ by deleting $v$ and its children; that is, $T' = T - V(T_v)$ where recall that $T_v$ denotes the maximal subtree of $T$ at $v$ induced by $D[v]$. Let $T'$ have order~$n'$, and so $n' = n - d_T(v) = n - \ell$. Since $\diam(T) \ge 4$, we note that $n' \ge 3$. Applying the inductive hypothesis to the tree $T'$, $F_t(T') \le ((\Delta' - 1)n' + 1)/\Delta' \le ((\Delta - 1)n' + 1)/\Delta$. Further, if $F_t(T') = ((\Delta - 1)n' + 1)/\Delta$, then $\Delta' = \Delta$ and $T' \in \cT_{\Delta}$. Let $S'$ be a minimum TF-set in $T'$, and so $|S'| = F_t(T')$. We note that every child of $v$ is a leaf. Let $S$ be the set obtained from $S'$ by adding to it $v$ and all children of $v$ different from $u$; that is, $S = S' \cup (D[v] \setminus \{u\})$. The set $S$ is a TF-set of $T$, implying that
\1
\[
\begin{array}{lcl}
F_t(T) \, \le \, |S| &  = & |S'| + \ell - 1 \1 \\
& = & F_t(T') + \ell - 1 \1 \\
& \le & \frac{1}{\Delta} ( (\Delta - 1)n' + 1) + \ell - 1 \2 \\
& = & \frac{1}{\Delta} ( (\Delta - 1)(n-\ell) + 1) + \ell - 1 \2  \\
& = & \frac{1}{\Delta} ( (\Delta - 1)n + 1)  +  \frac{1}{\Delta} (\ell - \Delta) \2  \\
& \le & \frac{1}{\Delta} ( (\Delta - 1)n + 1), \end{array}
\]
\2
\noindent
which establishes the desired upper bound of the theorem. Suppose that
\[
F_t(T) = \frac{1}{\Delta} ( (\Delta - 1)n + 1)
\]

\noindent (and still $\Delta \ge 3$ and $\diam(T) \ge 4$). Then, we must have equality throughout the above inequality chain, implying that $F_t(T') = \frac{1}{\Delta} ( (\Delta - 1)n' + 1)$ and $d_T(v) = \ell = \Delta$. By the inductive hypothesis, $T' \in \cT_{\Delta}$. If the parent $w$ of $v$ in $T$ is a central vertex of one of the underlying trees of $T'$, then it would have degree~$\Delta + 1$ in $T$, a contradiction. Hence, $w$ is a leaf in one of the underlying trees of $T' \in \cT_{\Delta}$. Let $T_1,\ldots,T_k$ denote the underlying trees of $T'$, and let $v_i$ be the central vertex of the tree $T_i$ for $i \in [k]$. Further, let $w$ belong to the subtree $T_j$, where $j \in [k]$. As observed earlier, $w$ is a leaf of the underlying tree $T_j$ of $T'$.

Since $T' \in \cT_{\Delta}$, we note that $T_1 \cong K_{1,\Delta}$, and if $k \ge 2$, then $T_i \cong K_{1,\Delta - 1}$ for $i \in [k] \setminus \{1\}$. Further, for $i \in [k]$, the central vertex $v_i$ of the star $T_i$ is a strong support vertex of degree~$\Delta$ in the tree $T'$, and the set $\{v_1,\ldots,v_k\}$ is an independent set in $T'$.

Suppose that the central vertex $v_j$ of $T_j$ is not a strong support vertex in $T$. Since $v_j$ is a strong support vertex in $T'$, this implies that $v_j$ has precisely two leaf neighbors in $T'$, one of which is necessarily the vertex $w$. Thus, $w$ is a leaf in $T'$ but has degree~$2$ in $T$, with $v$ and $x$ as its neighbors. This in turn implies that for $i \in [k] \setminus \{j\}$, every leaf neighbor of the vertex $v_i$ in $T'$ is also a leaf neighbor of $v_i$ in $T$, and so the vertex $v_i$ is a strong support vertex of degree~$\Delta$ in the tree $T$. For $i \in [k]$, let $v_i'$ be a leaf neighbor of $v_i$ in the tree $T$. We note that $v_j'$ is the unique leaf neighbor of $v_j$ in $T$. Let $L = \cup_{i=1}^k \{v_i'\}$.

We now consider the set $S = V(T) \setminus (L \cup \{u,x\})$. We note that the vertex $v$ and all its neighbors different from $u$ belong to $S$. Thus, playing the vertex $v$ in the first step of the forcing process starting with the initial set $S$ of colored vertices, the vertex $v$ forces its leaf neighbor $u$ to be colored. In the second step of the forcing process, we play the vertex $w$ which forces the vertex $x = v_j$ to be colored. At this stage of the forcing process, we note that all vertices of $V(T)$ are colored, except for the $k$ leaves $v_1',\ldots,v_k'$. Since $\{v_1,\ldots,v_k\}$ is an independent set in $T$, we now simply play each of the vertices $v_1, \ldots, v_k$ in turn in the forcing process, thereby forcing the leaves $v_1',\ldots,v_k'$ to be colored. In this way, we color all of $V(T)$ starting with the initial set $S$. Further, since $T[S]$ contains no isolated vertex, the set $S$ is therefore a TF-set. We note that $S$ contains all, except for $k + 1$, vertices of $V(T')$. By Lemma~\ref{l:lem1}(b), $F_t(T') = n' - k$, implying that $|S \cap V(T')| = n' - (k+1) < F_t(T')$. Thus,
\1
\[
\begin{array}{lcl}
F_t(T) \, \le \, |S| & = & |S \cap V(T')| + \Delta - 1 \1 \\
& < & F_t(T') + \Delta - 1 \1 \\
& = & \frac{1}{\Delta} ( (\Delta - 1)n' + 1) + \Delta - 1 \2 \\
& = & \frac{1}{\Delta} ( (\Delta - 1)(n - \Delta) + 1) + \Delta - 1 \2  \\
& = & \frac{1}{\Delta} ( (\Delta - 1)n + 1), \end{array}
\]

\noindent
a contradiction. Hence, the central vertex $v_j$ of $T_j$ is a strong support vertex in $T$. We now let $T_{k+1} = T_v$, where as defined earlier $T_v$ is the maximal subtree of $T$ at $v$ induced by $D[v]$, and we let $v_{k+1} = v$. We note that $T_{k+1} \cong K_{1,\Delta - 1}$. Further, we note that the central vertex $v_i$ of the star $T_i$ is a strong support vertex of degree~$\Delta$ in the tree $T$ for all $i \in [k+1]$, and the set $\{v_1,v_2,\ldots,v_{k+1}\}$ is an independent set in $T$. Thus, $T \in \cT_{\Delta}$. Conversely, if $T \in \cT_{\Delta}$, then by Lemma~\ref{l:lem1}, $F_t(T) = \frac{1}{\Delta} ( (\Delta - 1)n + 1)$. This completes the proof of Theorem~\ref{t:tree1}.~\qed

\section{Proof of Theorem~\ref{t:tree0}}
\label{S:tree0}

For all graphs of order~$n$ with maximum degree~$\Delta$, we note that $n \ge \Delta + 1$, implying that
\begin{equation}
\label{Eq1}
\frac{(\Delta - 1)n + 1}{\Delta} \le \left( \frac{\Delta}{\Delta + 1} \right) n.
\end{equation}

\noindent
Further, equality holds in Inequality~(\ref{Eq1}) if and only if $n = \Delta + 1$. Thus, the upper bound of Theorem~\ref{t:tree0} follows as an immediate consequence of the upper bound of Theorem~\ref{t:tree1}. Moreover, if $T$ is a tree for which equality holds in Inequality~(\ref{Eq1}), then $n = \Delta + 1$, implying that $T = K_{1,\Delta}$.

\section{Proof of Theorem~\ref{t:tree2}}
\label{S:cF}

In this section, we present a proof of Theorem~\ref{t:tree2}. For this purpose, we present a series of preliminary lemmas which will be used in our subsequent argument to establish the desired lower bound.

\begin{lem}\label{lemPaths}
Let $G$ be an isolate-free graph that contains an edge $e$ incident with a vertex of degree at most~$2$. If $G'$ is obtained from $G$ by subdividing the edge $e$ any number of times, then $F_t(G) = F_t(G')$.
\end{lem}
\proof Let $G$ be an isolate-free graph that contains an edge $e = uv$, where $v$ has degree~$1$ or~$2$ in $G$. Let $G'$ be obtained from $G$ by subdividing the edge $e$ any number of times. We may assume the edge $e$ is subdivided at least once, for otherwise $G' = G$ and the result is immediate. Let $P \colon u u_1 \ldots u_k v$ denote the resulting $(u,v)$-path in $G'$, where $k \ge 1$. If $d_G(v) = 2$, then let $w$ denote its neighbor different from $u$ in~$G$.

We first show that $F_t(G') \le F_t(G)$. Let $S \subseteq V(G)$ be a minimum TF-set of $G$, and so $|S| = F_t(G)$. Suppose that $v \in S$. Since $S$ is a TF-set of $G$, the graph $G[S]$ contains no isolated vertex. In particular, $S$ contains a neighbor of $v$. Suppose that both $u$ and $w$ belong to $S$. If $v$ is the only neighbor of $u$ in $S$, then the set $S \setminus \{u\}$ is a TF-set of $G$, contradicting the minimality of $S$. Hence, $u$ has at least two neighbors in $S$. Analogously, $w$ has at least two neighbors in $S$. Thus, the set $(S \setminus \{v\}) \cup \{u_1\}$ is a TF-set of $G'$, and so $F_t(G') \le |S| = F_t(G)$, as desired. If $u \in S$ and $w \notin S$, then once again the set $(S \setminus \{v\}) \cup \{u_1\}$ is a TF-set of $G'$, and $F_t(G') \le F_t(G)$, as desired. If $u \notin S$ and $w \in S$, then the set $S$ is a TF-set of $G'$, and so $F_t(G') \le |S| = F_t(G)$, as desired. Hence we may assume that $v \notin S$.

Since $S$ is a TF-set of $G$, there is a sequence $s \colon x_1,\ldots,x_t$ of played vertices in the forcing process that results in all $V(G)$ colored, where $x_i$ denotes the forcing colored vertex played in the $i$th step of the process. In particular, $v = x_\ell$ for some integer $\ell$ where $\ell \in [t]$. Before the vertex $v$ is colored, at least one neighbor of $v$ is colored. If $u$ is already colored before $v$, then starting with the same initial colored set $S$, the sequence obtained from $s$ by replacing the vertex $x_\ell$ by the subsequence $u_1, \ldots, u_k, v$ results in a sequence of played vertices in the forcing process that results in all $V(G')$ colored. If $u$ is colored after $v$, then $w$ is necessarily colored before $v$. In this case,  starting with the same initial colored set $S$, the sequence obtained from $s$ by replacing the vertex $x_\ell$ by the subsequence $v, u_k, u_{k-1}, \ldots, u_1$ results in a sequence of played vertices in the forcing process that results in all $V(G')$ colored. Thus, once again $F_t(G') \le |S| = F_t(G)$, as desired.

We first next that $F_t(G) \le F_t(G')$. Let $S' \subseteq V(G)$ be a minimum TF-set of $G'$, and so $|S| = F_t(G')$. If $|S' \cap V(P)| \le 1$, then since the graph $G'[S]$ contains no isolated vertex, the set $S'$ contains no vertex of $P$, except possibly for one of the ends of $P$, namely the vertex $u$ or the vertex $v$.
In both cases, the set $S'$ is a TF-set of $G$, and so $F_t(G) \le |S'| = F_t(G')$, as desired. If $|S' \cap V(P)| \ge 2$, then the set $(S' \setminus V(P)) \cup \{u,v\}$ is a TF-set of $G$, and so $F_t(G) \le |S'| - |S' \cap V(P)| + 2 \le |S'| = F_t(G')$, as desired.

Thus, $F_t(G') \le F_t(G)$ and $F_t(G) \le F_t(G')$. Consequently, $F_t(G) = F_t(G')$.~\qed

\begin{figure}[htb]
\begin{center}
\begin{tikzpicture}[scale=.8,style=thick,x=1cm,y=1cm]
\def\vr{2.5pt} 
%
\path (4,1) coordinate (u1);
\path (4.5,2) coordinate (v1);
\path (5,1) coordinate (u2);
\path (5.5,2) coordinate (v);
\path (5.15,2) coordinate (z1);
\path (5.7,2) coordinate (z2);
\path (6.25,2) coordinate (z3);
\path (6.8,2) coordinate (z4);
\path (5.5,1) coordinate (u9);
\path (7,1) coordinate (u3);
\path (7,-1) coordinate (u8);
\path (7.5,2) coordinate (v2);
\path (8,1) coordinate (u4);
\path (7,0) coordinate (x4);
\path (8,-1) coordinate (u7);
\path (8,0) coordinate (x3);
\path (5,0) coordinate (u5);
\path (5,-1) coordinate (x2);
\path (4, 0) coordinate (u6);
\path (4,-1) coordinate (x1);
\draw (v1) -- (v2);
\draw (v1) -- (u1);
\draw (v1) -- (u2);
\draw (v2) -- (u3);
\draw (v2) -- (u4);
\draw (u2) -- (u5);
\draw (u1) -- (u6);
\draw (u4) -- (u7);
\draw (u3) -- (u8);
\draw (x1) -- (u6);
\draw (x2) -- (u5);
\draw (x3) -- (u7);
\draw (x4) -- (u3);
\draw (v1) [fill=black] circle (\vr);
\draw (v2) [fill=black] circle (\vr);
\draw (z1) [fill=black] circle (\vr);
\draw (z2) [fill=black] circle (\vr);
\draw (z3) [fill=black] circle (\vr);
\draw (z4) [fill=black] circle (\vr);
\draw (u1) [fill=black] circle (\vr);
\draw (u2) [fill=black] circle (\vr);
\draw (u3) [fill=black] circle (\vr);
\draw (u4) [fill=black] circle (\vr);
\draw (u5) [fill=black] circle (\vr);
\draw (u6) [fill = black] circle (\vr);
\draw (u7) [fill = black] circle (\vr);
\draw (u8) [fill = black] circle (\vr);
\draw (x1) [fill = black] circle (\vr);
\draw (x2) [fill = black] circle (\vr);
\draw (x3) [fill = black] circle (\vr);
\draw (x4) [fill = black] circle (\vr);
\draw (6.5,-2) node {(a) $T$};
\path (12.5,0) coordinate (u1);
\path (13.5,0) coordinate (u2);
\path (14.5,0) coordinate (u5);
\path (15.5,0) coordinate (u6);

\path (13,1) coordinate (v1);
\path (14,1) coordinate (v2);
\path (15,1) coordinate (v3);

\draw (v1) -- (u1);
\draw (v1) -- (u2);
\draw (v3) -- (u5);
\draw (v3) -- (u6);
\draw (v1) -- (v3);

\draw (v1) [fill=black] circle (\vr);
\draw (v2) [fill=black] circle (\vr);
\draw (v3) [fill=black] circle (\vr);
\draw (u1) [fill=black] circle (\vr);
\draw (u2) [fill=black] circle (\vr);
\draw (u5) [fill=black] circle (\vr);
\draw (u6) [fill=black] circle (\vr);
\draw (14,-2) node {(b) $\text{trim}(T)$};
\end{tikzpicture}
\end{center}
\vskip -0.7 cm
\caption{A tree $T$ and its trimmed tree $\trim(T)$} \label{fig2}
\end{figure}
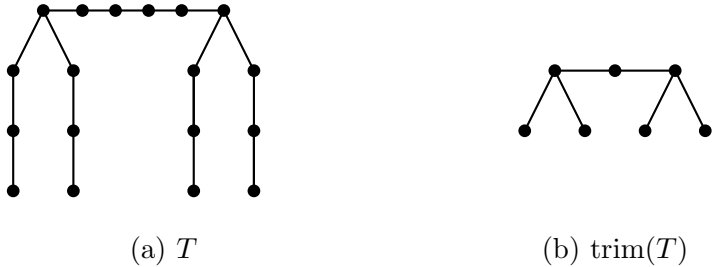

\medskip
The contraction of an edge $e = xy$ in a graph $G$ is the graph obtained from $G$ by replacing the vertices $x$ and $y$ by a new vertex and joining this new vertex to all vertices that were adjacent to $x$ or $y$ in $G$. Given a non-trivial tree $T$, the \emph{trimmed tree} of $T$, denoted $\trim(T)$, is the tree obtained from $T$ by iteratively contracting edges with one of its incident vertices of degree exactly~$2$ and with the other incident vertex of degree at most~$2$ until no such edge remains. We note that if the original tree $T$ is a path, then $\trim(T)$ is a path $P_2$, while if $T$ is not a path, then every edge in $\trim(T)$ is incident with a vertex of degree at least~$3$. In particular, if $T$ is not a path, then every support vertex in $\trim(T)$ has degree at least~$3$. As an illustration, the trimmed tree $\trim(T)$ of the tree $T$ shown in Figure~\ref{fig2}(a) is shown in Figure~\ref{fig2}(b).

\bigskip

Since every non-trivial tree $T$ can be reconstructed from its trimmed tree $\trim(T)$ by applying a sequence of subdivisions of edges incident with a vertex of degree at most~$2$, as an immediate consequence of Lemma~\ref{lemPaths} we note that $F_t(T) = F_t(\trim(T))$. We remark that the number of leaves in $T$ is equal to the number of leaves in $\trim(T)$. We state these properties of a trimmed tree formally as follows.

\begin{lem} \label{lemTrees}
If $T$ is a non-trivial tree, then the following hold. \\
\indent {\rm (a)} $F_t(T) = F_t(\trim(T))$. \\
\indent {\rm (b)} The trees $T$ and $\trim(T)$ have the same number of leaves.
\end{lem}

We proceed further by constructing a family $\cF$ of trees with small total forcing number.

\begin{figure}[htb]

\begin{center}
\begin{tikzpicture}[scale=.8,style=thick,x=1cm,y=1cm]
\def\vr{2.5pt} 
\path (1.7,1.15) coordinate (v);
\path (1.7,.15) coordinate (w);
%
\draw (v) -- (w);
\draw (v) [fill=black] circle (\vr);
\draw (w) [fill=black] circle (\vr);
\draw [rounded corners] (0,-.2) rectangle (2,1.5);
\draw (.4,0.75) node {$T'$};
\draw (-1,.75) node {$\cO_1$:};
\draw (-2,.75) node {(a)};
\draw[anchor = east] (v) node {$u$};
\draw[anchor = east] (w) node {$v$};
\draw (3,0.75) node {$\mapsto$};
\path (5.7,1.15) coordinate (v);
\path (6.7,0.65) coordinate (x);
\path (5.7,0.15) coordinate (w);
%
\draw (v) -- (x);
\draw (w) -- (x);
\draw (v) [fill=black] circle (\vr);
\draw (x) [fill=white] circle (\vr);
\draw (w) [fill=black] circle (\vr);
\draw [rounded corners] (4,-.2) rectangle (6,1.5);
\draw[anchor = east] (v) node {$u$};
\draw[anchor = west] (x) node {$w$};
\draw[anchor = east] (w) node {$v$};
\end{tikzpicture}
\end{center}

\vskip 0.2 cm

\begin{center}
\begin{tikzpicture}[scale=.8,style=thick,x=1cm,y=1cm]
\def\vr{2.5pt} 
\path (1.1,1.05) coordinate (u);
\path (1.7,1.05) coordinate (w);
\path (1.4,.25) coordinate (v);
\path (0.8,.15) coordinate (a);
\path (0.8,.65) coordinate (b);
%
\draw (v) -- (u);
\draw (v) -- (w);
\draw (v) -- (a);
\draw (v) -- (b);
\draw (v) [fill=black] circle (\vr);
\draw (u) [fill=black] circle (\vr);
\draw (w) [fill=black] circle (\vr);
\draw [rounded corners] (0,-.2) rectangle (2.2,1.5);
\draw (.4,0.75) node {{\small $T'$}};
\draw (-1,.75) node {$\cO_2$:};
\draw (-2,.75) node {(b)};
\draw[anchor = north] (v) node {{\small $v$}};
\draw[anchor = south] (u) node {{\small $u$}};
\draw[anchor = south] (w) node {{\small $w$}};
\draw (3,0.75) node {$\mapsto$};
\path (5.1,1.05) coordinate (u);
\path (5.7,1.05) coordinate (w);
\path (5.4,.25) coordinate (v);
\path (4.8,.15) coordinate (a);
\path (4.8,.65) coordinate (b);
\path (6.75,.25) coordinate (x);
%
\draw (v) -- (u);
\draw (v) -- (w);
\draw (v) -- (a);
\draw (v) -- (b);
\draw (v) -- (x);
\draw (v) [fill=black] circle (\vr);
\draw (u) [fill=black] circle (\vr);
\draw (w) [fill=black] circle (\vr);
\draw (x) [fill=white] circle (\vr);
\draw [rounded corners] (4,-.2) rectangle (6,1.5);
\draw[anchor = north] (v) node {{\small $v$}};
\draw[anchor = south] (u) node {{\small $u$}};
\draw[anchor = south] (w) node {{\small $w$}};
\draw[anchor = west] (x) node {{\small $x$}};
\end{tikzpicture}
\end{center}

\vskip 0.2 cm

\begin{center}
\begin{tikzpicture}[scale=.8,style=thick,x=1cm,y=1cm]
\def\vr{2.5pt} 
\path (1.1,0.85) coordinate (u);
\path (1.7,0.65) coordinate (w);
\path (1.4,.25) coordinate (v);
\path (0.8,.15) coordinate (a);
\path (0.8,.65) coordinate (b);
%
\draw (v) -- (u);
\draw (v) -- (w);
\draw (v) -- (a);
\draw (v) -- (b);
\draw (v) [fill=black] circle (\vr);
\draw (u) [fill=black] circle (\vr);
\draw (w) [fill=black] circle (\vr);
\draw [rounded corners] (0,-.2) rectangle (2.2,1.5);
\draw (.4,0.75) node {{\small $T'$}};
\draw (-1,.75) node {$\cO_3$:};
\draw (-2,.75) node {(c)};
\draw[anchor = north] (v) node {{\small $w$}};
\draw[anchor = south] (u) node {{\small $u$}};
\draw[anchor = south] (w) node {{\small $v$}};
\draw (3,0.75) node {$\mapsto$};
\path (5.1,0.85) coordinate (u);
\path (5.7,0.65) coordinate (w);
\path (5.4,.25) coordinate (v);
\path (4.8,.15) coordinate (a);
\path (4.8,.65) coordinate (b);
\path (6.75,.25) coordinate (x);
\path (6.75,1.05) coordinate (y);
%
\draw (v) -- (u);
\draw (v) -- (w);
\draw (v) -- (a);
\draw (v) -- (b);
\draw (w) -- (x);
\draw (w) -- (y);
\draw (v) [fill=black] circle (\vr);
\draw (u) [fill=black] circle (\vr);
\draw (w) [fill=black] circle (\vr);
\draw (x) [fill=white] circle (\vr);
\draw (y) [fill=white] circle (\vr);
\draw [rounded corners] (4,-.2) rectangle (6,1.5);
\draw[anchor = north] (v) node {{\small $w$}};
\draw[anchor = south] (u) node {{\small $u$}};
\draw[anchor = south] (w) node {{\small $v$}};
\draw[anchor = west] (x) node {{\small $x$}};
\draw[anchor = west] (y) node {{\small $y$}};
\end{tikzpicture}
\end{center}

\vskip 0.2 cm

\begin{center}
\begin{tikzpicture}[scale=.8,style=thick,x=1cm,y=1cm]
\def\vr{2.5pt} 
\hspace*{0.25cm}
\path (1.6,0.65) coordinate (v);
\path (0.8,.25) coordinate (a);
\path (0.8,.95) coordinate (c);
%
\draw (v) -- (a);
\draw (v) -- (c);
\draw (v) [fill=black] circle (\vr);
\draw [rounded corners] (0,-.2) rectangle (2.2,1.5);
\draw (.4,0.75) node {{\small $T'$}};
\draw (-1,.75) node {$\cO_4$:};
\draw (-2,.75) node {(d)};
\draw[anchor = north] (v) node {{\small $v$}};
\draw (3,0.75) node {$\mapsto$};
\path (5.4,.65) coordinate (v);
\path (4.8,.25) coordinate (a);
\path (4.8,.95) coordinate (b);
\path (6.75,.65) coordinate (y);
\path (7.35,1.05) coordinate (x);
\path (7.35,.25) coordinate (z);
%
\draw (v) -- (a);
\draw (v) -- (b);
\draw (v) -- (y);
\draw (y) -- (x);
\draw (y) -- (z);
\draw (v) [fill=black] circle (\vr);
\draw (z) [fill=white] circle (\vr);
\draw (x) [fill=white] circle (\vr);
\draw (y) [fill=white] circle (\vr);
\draw [rounded corners] (4,-.2) rectangle (6,1.5);
\draw[anchor = north] (v) node {{\small $v$}};
\draw[anchor = west] (z) node {{\small $z$}};
\draw[anchor = west] (x) node {{\small $x$}};
\draw[anchor = north] (y) node {{\small $y$}};
\end{tikzpicture}
\end{center}

\vskip 0.2 cm

\begin{center}
\begin{tikzpicture}[scale=.8,style=thick,x=1cm,y=1cm]
\def\vr{2.5pt} 
\hspace*{0.75cm}
\path (1.6,0.65) coordinate (v);
\path (0.8,.25) coordinate (a);
\path (0.8,.95) coordinate (c);
%
\draw (v) -- (a);
\draw (v) -- (c);
\draw (v) [fill=black] circle (\vr);
\draw [rounded corners] (0,-.2) rectangle (2.2,1.5);
\draw (.4,0.75) node {{\small $T'$}};
\draw (-1,.75) node {$\cO_5$:};
\draw (-2,.75) node {(e)};
\draw[anchor = north] (v) node {{\small $v$}};
\draw (3,0.75) node {$\mapsto$};
\path (5.4,.65) coordinate (v);
\path (4.8,.25) coordinate (a);
\path (4.8,.95) coordinate (b);
\path (6.75,.65) coordinate (v1);
\path (6.75,.15) coordinate (u1);
\path (7.45,.65) coordinate (v2);
\path (8.15,.65) coordinate (v3);
\path (8.15,.15) coordinate (u3);
\path (8.85,.65) coordinate (v4);
%
\draw (v) -- (a);
\draw (v) -- (b);
\draw (v) -- (v1);
\draw (v1) -- (v2);
\draw (v2) -- (v3);
\draw (v3) -- (v4);
\draw (v1) -- (u1);
\draw (v3) -- (u3);
\draw (v) [fill=black] circle (\vr);
\draw (v1) [fill=white] circle (\vr);
\draw (v2) [fill=white] circle (\vr);
\draw (v3) [fill=white] circle (\vr);
\draw (v4) [fill=white] circle (\vr);
\draw (u1) [fill=white] circle (\vr);
\draw (u3) [fill=white] circle (\vr);
\draw [rounded corners] (4,-.2) rectangle (6,1.5);
\draw[anchor = north] (v) node {{\small $v$}};
\draw[anchor = south] (v1) node {{\small $v_1$}};
\draw[anchor = south] (v2) node {{\small $v_2$}};
\draw[anchor = south] (v3) node {{\small $v_3$}};
\draw[anchor = south] (v4) node {{\small $v_4$}};
\draw[anchor = north] (u1) node {{\small $u_1$}};
\draw[anchor = north] (u3) node {{\small $u_3$}};
\end{tikzpicture}
\end{center}

\vskip -0.4cm
\caption{The operations $\cO_1$, $\cO_2$, $\cO_3$, $\cO_4$, and $\cO_5$.} \label{f:cO1O5}
\end{figure}
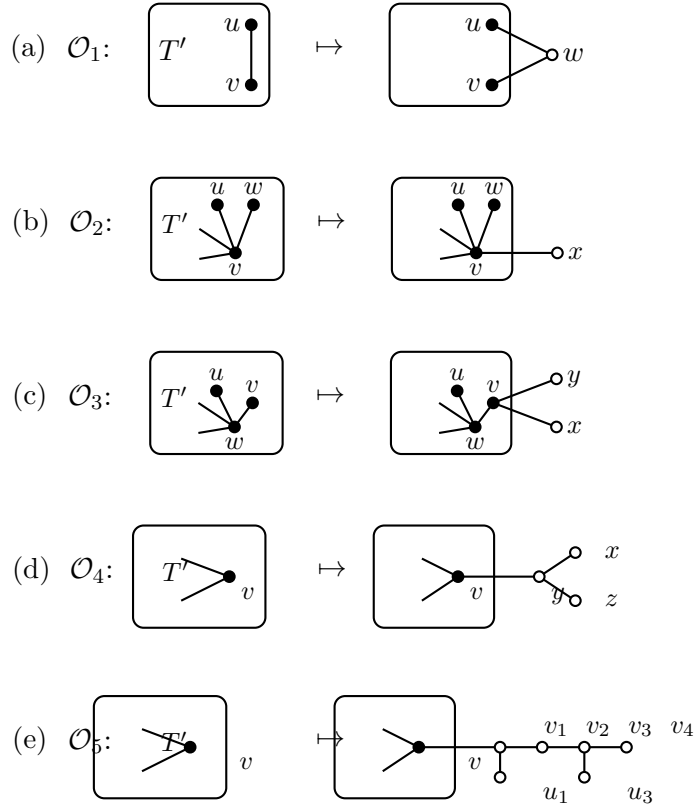

\noindent\textbf{The family $\cF$.}
Let $\cF$ be the family of trees that contains a path $P_2$ and is closed under the five operations $\cO_1, \cO_2, \ldots, \cO_5$ below, which extend a tree $T'$ to a new tree~$T$. In Fig.~\ref{f:cO1O5}, the vertices of $T'$ are colored black and the new vertices of $T$ are colored white.

\noindent Operation $\cO_1$: If $uv$ is an edge of $T'$ where at least one of $u$ and $v$ has degree at most~$2$ in $T'$, then $T$ is obtained from $T'$ by subdividing the edge $uv$ once. See Fig~\ref{f:cO1O5}(a), where $w$ denotes the new vertex (of degree~$2$ in $T$) obtained from subdividing the edge $uv$.

\noindent Operation $\cO_2$: If $v$ is a strong support vertex in $T'$, then $T$ is obtained from $T'$ by adding an additional pendant edge to the vertex $v$. See Fig~\ref{f:cO1O5}(b), where $u$ and $w$ are leaf neighbors of $v$ in $T'$ and $x$ is the new leaf added to $T'$.

\noindent Operation $\cO_3$: If $w$ is a strong support vertex in $T'$ and $v$ is a leaf neighbor of $w$, then $T$ is obtained from $T'$ by adding two pendant edges to $v$. See Fig~\ref{f:cO1O5}(c), where $u$ and $v$ are leaf neighbors of $w$ in $T'$, and $vx$ and $vy$ are the pendant edges added to $v$.

\noindent Operation $\cO_4$: If $v$ is a vertex of degree at least~$2$ in $T'$, then $T$ is obtained from $T'$ by adding a path $P_3$ and joining $v$  to the central vertex of the path. See Fig~\ref{f:cO1O5}(d) where $xyz$ is the added path and $vy$ the added edge.

\noindent Operation $\cO_5$: If $v$ is a vertex of degree at least~$2$ in $T'$, then $T$ is obtained from $T'$ by adding a star $K_{1,3}$ with one edge subdivided twice and adding an edge joining the resulting support vertex of degree~$2$ to the vertex $v$. See Fig~\ref{f:cO1O5}(e) where the added star has central vertex $v_3$ with leaf neighbors $u_1$, $u_3$ and $v_4$, and where the edge $u_1v_3$ of the star is subdivided twice resulting in the path $u_1v_1v_2v_3$ with $v_1$ a support vertex of degree~$2$ in the subdivided star.

In operations $\cO_2$, $\cO_3$, $\cO_4$  and $\cO_5$ illustrated in Figure~\ref{f:cO1O5}, we call the vertex $v$ the \emph{link vertex} of $T'$. We are now in a position to establish the following lower bound on the total forcing number of a tree in terms of the number of leaves in the tree.

\begin{lem}
\label{l:Leaves1}
If $T$ is a non-trivial tree with $n_1$ leaves, then $F_t(T) \ge n_1$. Further, if $F_t(T) = n_1$, then $T \in \cF$.
\end{lem}
\proof  We proceed by induction on the order~$n \ge 2$ of a tree $T$ with $n_1$ leaves. If $n = 2$, then $T = P_2 \in \cF$ and $F_t(T) = 2 = n_1$. This establishes the base case. Let $n \ge 3$ and assume that if $T'$ is a non-trivial tree of order~$n'$ where $n' < n$ having $n_1'$ leaves, then $F_t(T') \ge n_1'$ and that if $F_t(T') = n_1'$, then $T' \in \cF$. Let $T$ be a tree of order~$n$ with $n_1$ leaves. We proceed further with the following series of claims.

\begin{claim}
\label{c:claim1}
If $T \ne \trim(T)$, then $F_t(T) \ge n_1$ and if $F_t(T) = n_1$, then $T \in \cF$.
\end{claim}
\proof Suppose that the trimmed tree, $\trim(T)$, is different from $T$. This implies that $T$ contains an edge $e$ with one of its incident vertices of degree exactly~$2$ and with the other incident vertex of degree at most~$2$. Let $T'$ be obtained from $T$ by contracting the edge $e$. By Lemma~\ref{lemPaths}, $F_t(T) = F_t(T')$. Let $T'$ have $n_1'$ leaves. We note that $T'$ has order~$n-1$, and that $n_1 = n_1'$. Applying the inductive hypothesis to the tree~$T'$, $F_t(T) = F_t(T') \ge n_1' = n_1$. Further, if $F_t(T) = n_1$, then $F_t(T') = n_1'$ and by the inductive hypothesis, $T' \in \cF$. In this case, we can restore the tree $T$ by applying operation~$\cO_1$ to the tree $T'$, implying that $T \in \cF$.~\smallqed

\medskip
By Claim~\ref{c:claim1}, we may assume that $T = \trim(T)$, for otherwise the desired result follows. With this assumption, we note that every edge in $T$ is incident with a vertex of degree at least~$3$. In particular, every support vertex in $T$ has degree at least~$3$.

\begin{claim}
\label{c:claim2}
If $T$ contains a support vertex with at least three leaf neighbors, then $F_t(T) \ge n_1$ and if $F_t(T) = n_1$, then $T \in \cF$.
\end{claim}
\proof Suppose that $T$ contains a support vertex $v$ with at least three leaf neighbors. Let $S$ be a minimum TF-set of $F$, and so $|S| = F_t(F)$. By Observation~\ref{ob:2}, the set $S$ contains the vertex $v$ and all except possibly one leaf neighbor of~$v$. Let $u$ and $u'$ be two distinct leaf neighbors of $v$. If $S$ contains every leaf neighbor of~$v$, then by the minimality of $S$, there is a neighbor, $w$ say, of $v$ not in $S$. In this case, replacing the vertex $u$ in $S$ with the vertex $w$ produces a new minimum TF-set of $F$. Hence, renaming the leaf neighbors of $v$ if necessary, we may assume that $u \notin S$. We now consider the tree $T' = T - u'$. Let $T'$ have $n_1'$ leaves, and so $n_1' = n_1 - 1$. The set $S \setminus \{u'\}$ is necessarily a TF-set of $T'$, implying that $F_t(T') \le |S| - 1 = F_t(T) - 1$. Applying the inductive hypothesis to the tree $T'$, we therefore have $n_1 - 1 = n_1' \le F_t(T') \le F_t(T) - 1$, implying that $F_t(T) \ge n_1$. Further, suppose that $F_t(T) = n_1$. In this case, $F_t(T') = n_1'$. Applying the inductive hypothesis to the tree $T'$, we have $T' \in \cF$. We note that in the tree $T'$ the vertex $v$ is a strong support vertex. Hence, $T$ can be obtained from the tree $T'$ by applying operation~$\cO_2$ with $v$ as the link vertex, and so $T \in \cF$.~\smallqed

\medskip
By Claim~\ref{c:claim2}, we may assume that every support vertex in $T$ and has at most two leaf neighbors. Recall that every support vertex in $T$ has degree at least~$3$. With these assumptions, we note that $T$ is not a star, and so $\diam(T) \ge 3$.

\begin{claim}
\label{c:claim3}
If $\diam(T) = 3$, then $F_t(T) = n_1$ and $T \in \cF$.
\end{claim}
\proof Suppose that $\diam(T) = 3$, and so $T \cong S(r,s)$ is a double star, where $1 \le r \le s$. Since  $T = \trim(T)$, both vertices of $T$ that are not leaves are strong support vertices, and so $r \ge 2$. By assumption, every support vertex in $T$ has at most two leaf neighbors, and so $s \le 2$. Thus, $r = s = 2$ and $T$ is the double star $S(2,2)$. By Observation~\ref{ob:2}, $F_t(T) \ge 4 = n_1$. If $\ell_1$ and $\ell_2$ are arbitrary leaves at distance~$3$ apart in $T$, then the set $V(T) \setminus \{\ell_1,\ell_2\}$ is a TF-set of $T$, and so $F_t(T) \le |V(T)| - 2 = 4 = n_1$. Consequently, $F_t(T) = n_1$. Further, we note that the double star $T$ can be constructed from a path $P_2$ by first applying operation $\cO_1$ and then applying operation $\cO_4$. Thus, $T \in \cF$.~\smallqed

\medskip
By Claim~\ref{c:claim3}, we may assume that $\diam(T) \ge 4$, for otherwise the desired result holds. Let $u$ and $r$ be two vertices at maximum distance apart in $T$. Necessarily, $u$ and $r$ are leaves and $d(u,r) = \diam(T)$.
We now root the tree $T$ at the vertex $r$. Let $v$ be the parent of $u$, let $w$ be the parent of $v$, let $x$ be the parent of $w$, and let $y$ be the parent of $x$. Possibly, $y = r$. We note that every child of $v$ is a leaf. Since $T = \trim(T)$, we note that $d_T(v) \ge 3$. By assumption, every support vertex in $T$ has at most two leaf neighbors. Thus, $d_T(v) \le 3$. Consequently, $d_T(v) = 3$. Let $u'$ be the child of $v$ different from $u$. Thus, $N_T(v) = \{u,u',w\}$, where recall that $w$ is the parent of $v$ in $T$.

Let $S$ be a minimum TF-set of $F$, and so $|S| = F_t(F)$. If $S$ contains both children $u$ and $u'$ of $v$, then by the minimality of $S$, the vertex $w \notin S$. In this case, replacing $u$ in $S$ with the vertex $w$ produces a new minimum TF-set of $F$. Hence, renaming the children of $v$ if necessary, we may assume that $u \notin S$. By Observation~\ref{ob:2}, the set $S$ therefore contains the vertex $v$ and its child $u'$. Thus, $\{u',v\} \subseteq S$ and $u \notin S$.

\begin{claim}
\label{c:claim5}
If $d_T(w) \ge 3$, then $F_t(T) \ge n_1$ and if $F_t(T) = n_1$, then $T \in \cF$.
\end{claim}
\proof Suppose that $d_T(w) \ge 3$. We now consider the tree $T' = T - \{u,u',v\}$. Let $T'$ have $n_1'$ leaves, and so $n_1' = n_1 - 2$. Let $S' = S \setminus \{u',v\}$. If $S'$ is a TF-set of $T'$, then, by the inductive hypothesis, $n_1 - 2 = n_1' \le F_t(T') \le |S'| = |S| - 2 = F_t(T) - 2$, implying that $F_t(T) \ge n_1$.  Further, suppose that $F_t(T) = n_1$. In this case, $F_t(T') = n_1'$. Applying the inductive hypothesis to the tree $T'$, we have $T' \in \cF$. We note that in the tree $T'$ the vertex $w$ has degree at least~$2$, and therefore $T$ can be obtained from the tree $T'$ by applying operation~$\cO_4$ with $w$ as the link vertex, and so $T \in \cF$.

Hence, we may assume that the set $S'$ is not a TF-set of $T'$, implying that $S$ contains $w$ but no neighbor of $w$ except for its child~$v$. In particular, $S$ contains no neighbor of $w$ in $T'$. If a child $v'$ of $w$ different from $v$ is not a leaf, then identical arguments as shown with the vertex $v$ show that $d_T(v') = 3$ and that $v' \in S$. Thus, $S$ contains a neighbor of $w$ different from $v$, a contradiction. Hence, every child of $w$ different from $v$ is a leaf. If $w$ has at least two leaf neighbors, then by Observation~\ref{ob:2}, the set $S$ contains $w$ and all except possibly one leaf neighbor of $w$, implying once again that $S$ contains a neighbor of $w$ different from $v$, a contradiction. Therefore, $d_T(w) = 3$ and the child, $v'$ say, of $w$ different from $v$ is a leaf.

We now consider the tree $T'' = T - \{u,u'\}$. Let $T''$ have $n_1''$ leaves, and so $n_1'' = n_1 - 1$ noting that the vertex $v$ is a leaf in $T''$ but not in $T$. Recall that both $v$ and $w$ belong to $S$. The set $S \setminus \{u'\}$ is a TF-set of $T''$, implying that $F_t(T'') \le |S| - 1 = F_t(T) - 1$. Applying the inductive hypothesis to the tree $T'$, we therefore have $n_1 - 1 = n_1'' \le F_t(T'') \le F_t(T) - 1$, implying that $F_t(T) \ge n_1$. Further, suppose that $F_t(T) = n_1$. In this case, $F_t(T'') = n_1''$. Applying the inductive hypothesis to the tree $T''$, we have $T'' \in \cF$. We note that in the tree $T''$ the vertex $w$ is a strong support vertex with $v$ as one of its leaf neighbors. Hence, $T$ can be obtained from the tree $T''$ by applying operation~$\cO_3$ with $v$ as the link vertex, and so $T \in \cF$.~\smallqed

\medskip
By Claim~\ref{c:claim5}, we may assume that $d_T(w) = 2$, for otherwise the desired result holds. Since $T = \trim(T)$, this implies that $d_T(x) \ge 3$.

\begin{claim}
\label{c:claim6}
If $w \notin S$ or if $\{w,x\} \subset S$, then $F_t(T) \ge n_1 + 1$.
\end{claim}
\proof Suppose that $w \notin S$ or $\{w,x\} \subset S$. In this case, we consider the tree $T' = T - \{u,u',v\}$. Let $T'$ have $n_1'$ leaves, and so $n_1' = n_1 - 1$, noting that the vertex $w$ is a leaf in $T'$ but not in $T$. Let $S' = S \setminus \{u',v\}$. Since $w \notin S$ or $\{w,x\} \subset S$, the set $S'$ is a TF-set of $T'$, and so, by the inductive hypothesis, $n_1 - 1 = n_1' \le F_t(T') \le |S'| = |S| - 2 = F_t(T) - 2$, implying that $F_t(T) \ge n_1 + 1$.~\smallqed

\medskip
By Claim~\ref{c:claim6}, we may assume that $w \in S$ and $x \notin S$, for otherwise the desired result holds. Thus, by our earlier assumptions, $S \cap \{u,u',v,w,x\} = \{u',v,w\}$. Recall that $d_T(x) \ge 3$.

\begin{claim}
\label{c:claim7}
If $d_T(v) \ge 4$ or if $d_T(x) = 3$ and the child of $x$ different from $w$ is not a leaf,
then $F_t(T) \ge n_1$ and if $F_t(T) = n_1$, then $T \in \cF$.
\end{claim}
\proof Suppose that $d_T(v) \ge 4$ or $d_T(x) = 3$ and the child of $x$ different from $w$ is not a leaf.
If $x$ has at least two leaf neighbors, then by our earlier assumptions, $x$ has exactly two leaf neighbors. In this case, by Observation~\ref{ob:2} the set $S$ contains the vertex $x$ (and at least one leaf neighbor of~$x$). This contradicts our assumption that $x \notin S$. Hence, at most one child of $x$ is a leaf. By assumption, there is a child $w_1$ of $x$ different from $w$ of degree at least~$2$. We note that either $w_1$ has a grandchild or every child of $w_1$ is a leaf.

Suppose firstly that $w_1$ has a grandchild, say $u_1$. Let $v_1$ be the parent of $u_1$. Using analogous arguments as before with the vertices $v$ and $w$, we may assume that $d_T(v_1) = 3$ and $d_T(w_1) = 2$, for otherwise the desired result follows. Let $u_2$ denote the child of $v_1$ different from $u_1$. Analogously as before, we may assume that $S \cap \{u_1,u_2,v_1,w_1\} = \{u_1,v_1,w_1\}$, for otherwise the desired result follows. We now consider the tree $T' = T - \{u,u',v,w\}$, and let $T'$ have $n_1'$ leaves. Thus, $n_1' = n_1 - 2$. The set $S \setminus \{u',v,w\}$ is a TF-set of $T'$, and so, by the inductive hypothesis, $n_1 - 2 = n_1' \le F_t(T') \le |S| - 3 = F_t(T) - 3$, implying that $F_t(T) \ge n_1 + 1$. Hence, we may assume that $w_1$ has no grandchild, for otherwise $F_t(T) \ge n_1 + 1$. Thus, every child of $w_1$ is a leaf.

Since $T = \trim(T)$, we note that $d_T(w_1) \ge 3$. By assumption, every support vertex in $T$ has at most two leaf neighbors. Thus, $d_T(w_1) \le 3$. Consequently, $d_T(w_1) = 3$. Let $v_1$ and $v_2$ denote the two children of $w_1$. By Observation~\ref{ob:2}, the set $S$ contains the vertex $w_1$ and at least one of $v_1$ and $v_2$. If $S$ contains both $v_1$ and $v_2$, then replacing $v_1$ in $S$ with the vertex $x$ produces a new minimum TF-set of $F$ that contains the vertex $x$, contradicting our earlier assumptions. Hence, renaming $v_1$ and $v_2$ if necessary, we may assume that $S \cap \{v_1,v_2,w_1\} = \{v_1,w_1\}$. We now consider the tree $T'' = T - \{v_1,v_2,w_1\}$. Let $T''$ have $n_1''$ leaves. Thus, $n_1'' = n_1 - 2$. The set $S \setminus \{v_1,w_1\}$ is a TF-set of $T''$, and so, by the inductive hypothesis, $n_1 - 2 = n_1'' \le F_t(T'') \le |S| - 2 = F_t(T) - 2$, implying that $F_t(T) \ge n_1$. Further, if $F_t(T) = n_1$, then $F_t(T'') = n_1''$ and by the inductive hypothesis, $T'' \in \cF$. In this case, we can restore the tree $T$ by applying operation~$\cO_4$ to the tree $T''$ with $x$ as the link vertex, implying that $T \in \cF$.~\smallqed

\medskip
By Claim~\ref{c:claim7}, we may assume that $d_T(x) = 3$ and the child, $w'$ say, of $x$ different from $w$ is a leaf, for otherwise the desired result holds. Recall that by our earlier assumptions, $x \notin S$, implying that $w' \notin S$.

\begin{claim}
\label{c:claim8}
If $d_T(y) = 2$, then $F_t(T) \ge n_1 + 1$.
\end{claim}
\proof Suppose that $d_T(y) = 2$. In this case, we consider the tree $T' = T - \{u,u',v,w,w',x\}$; that is, $T'$ is the tree obtained from $T$ by deleting $x$ and all its descendants. Let $T'$ have $n_1'$ leaves. Thus, $n_1' = n_1 - 2$, noting that $y$ is a leaf in $T'$ but is not a leaf in $T$. The set $S \setminus \{u',v,w\}$ is necessarily a TF-set of $T'$, and so, by the inductive hypothesis, $n_1 - 2 = n_1' \le F_t(T') \le |S| - 3 = F_t(T) - 3$, implying that $F_t(T) \ge n_1 + 1$.~\smallqed

\medskip
By Claim~\ref{c:claim8}, we may assume that  $d_T(y) \ge 3$, for otherwise the desired result holds. We now consider the tree $T' = T - \{u,u',v,w,w',x\}$; that is, $T'$ is the tree obtained from $T$ by deleting $x$ and all its descendants. Let $T'$ have $n_1'$ leaves. Thus, $n_1' = n_1 - 3$. The set $S \setminus \{u',v,w\}$ is necessarily a TF-set of $T'$, and so, by the inductive hypothesis, $n_1 - 3 = n_1' \le F_t(T') \le |S| - 3 = F_t(T) - 3$, implying that $F_t(T) \ge n_1$. Further, if $F_t(T) = n_1$, then $F_t(T') = n_1'$ and by the inductive hypothesis, $T' \in \cF$. In this case, we can restore the tree $T$ by applying operation~$\cO_5$ to the tree $T'$ with $y$ as the link vertex, implying that $T \in \cF$. This completes the proof of Lemma~\ref{l:Leaves1}.~\qed

\medskip
We show next that every tree in the family~$\cF$ has total forcing number equal to the number of leaves in the tree.

\begin{lem}
\label{l:Leaves2}
If $T$ is a tree in the family~$\cF$ with $n_1$ leaves, then $F_t(T) = n_1$.
\end{lem}
\proof We proceed by induction on the order~$n \ge 2$ of a tree $T$ in the family~$\cF$ with $n_1$ leaves. If $n = 2$, then $T = P_2$ and $F_t(T) = 2 = n_1$. This establishes the base case. Let $n \ge 3$ and assume that if $T' \in \cF$ is a tree of order~$n'$ where $n' < n$ and with $n_1'$ leaves, then $F_t(T') = n_1'$. Let $T$ be a tree of order~$n$ in the family~$\cF$ with $n_1$ leaves.
By definition of the family~$\cF$, there is a sequence $T_0,T_1, \ldots, T_k$ of trees where $T_0 = P_2$, $T_k = T$ and for $i \in [k]$, the tree $T_i$ can be obtained from the tree $T_{i-1}$ by one of the five operations $\cO_1, \cO_2, \ldots, \cO_5$. Let $T' = T_{k-1}$. Hence, $T' \in \cF$ and the tree $T'$ has order less than~$n$. Let $T'$ have $n_1'$ leaves. By the inductive hypothesis, $F_t(T) = n_1'$. Let $S'$ be a minimum TF-set of $T'$.

Suppose that $T$ is obtained from $T'$ by applying operation $\cO_1$. In this case, $n_1 = n_1'$ and, by Lemma~\ref{lemPaths}, $F_t(T) = F_t(T')$, implying that $F_t(T) = n_1$.

Suppose that $T$ is obtained from $T'$ by applying operation $\cO_2$. In this case, $n_1 = n_1' + 1$. Adopting the notation of Figure~\ref{f:cO1O5}(b), by Observation~\ref{ob:2} the set $S'$ contains the vertex $v$ and all except possibly one leaf neighbor. Thus, the set $S' \cup \{x\}$ is a TF-set of $T$, and so $F_t(T) \le |S'| + 1 = F_t(T') + 1 = n_1' + 1 = n_1$.

Suppose that $T$ is obtained from $T'$ by applying operation $\cO_3$. In this case, $n_1 = n_1' + 1$. Adopting the notation of Figure~\ref{f:cO1O5}(c), by Observation~\ref{ob:2} we may choose $S'$ so that $\{v,w\} \subseteq S$. Thus, the set $S' \cup \{x\}$ is a TF-set of $T$, and so $F_t(T) \le |S'| + 1 = F_t(T') + 1 = n_1' + 1 = n_1$.

Suppose that $T$ is obtained from $T'$ by applying operation $\cO_4$. In this case, $n_1 = n_1' + 2$. Adopting the notation of Figure~\ref{f:cO1O5}(d), every minimum TF-set of $T'$ can be extended to a TF-set of $T$ by adding to it $x$ and $y$, and so $F_t(T) \le F_t(T') + 2 = n_1' + 2 = n_1$.

Suppose that $T$ is obtained from $T'$ by applying operation $\cO_5$. In this case, $n_1 = n_1' + 3$. Adopting the notation of Figure~\ref{f:cO1O5}(e), every minimum TF-set of $T'$ can be extended to a TF-set of $T$ by adding to it the vertices $v_2$, $v_3$ and $v_4$, and so $F_t(T) \le F_t(T') + 3 = n_1' + 3 = n_1$. In all the above cases, $F_t(T) \le n_1$. By Lemma~\ref{l:Leaves1}, $F_t(T) \ge n_1$. Consequently, $F_t(T) = n_1$.~\qed

As an immediate consequence of Lemmas~\ref{l:Leaves1} and~\ref{l:Leaves2}, we have Theorem~\ref{t:tree2}. Recall its statement.

\noindent \textbf{Theorem~\ref{t:tree2}}. \emph{If $T$ is a non-trivial tree with $n_1$ leaves, then $F_t(T) \ge n_1$, with equality if and only if $T\in \cF$.}

\section{Proof of Theorem~\ref{t:tree3}}
\label{S:cS}

In this section we present a proof of Theorem~\ref{t:tree3}. We first observe that a simple adaptation of the proof given for Lemma~\ref{lemTrees} yields an analogous result on forcing. We state this formally with the following lemma.

\begin{lem}\label{forcingTreeLem}
If $T$ is a non-trivial tree, then $F(T) = F(\trim(T))$.
\end{lem}

We proceed further by defining a family $\cH$ of trees as follows.

\medskip
\noindent\textbf{The family $\cH$.} Let $n \ge 2$ be an integer and let $\cH$ be the family of all trees $T$ of order~$n$ such that $T \cong P_n$ or $\trim(T) \cong K_{1,n-1}$ for $n \ge 3$.

We are now in a position to prove Theorem~\ref{t:tree3}. Recall its statement.

\medskip
\noindent \textbf{Theorem~\ref{t:tree3}}. \emph{If $T$ is a non-trivial tree, then $F_t(T) \ge F(T)+1$, with equality if and only if $T \in \cH$.}

\proof Let $T$ be a non-trivial tree with $n_1$ leaves. By Theorem~\ref{t:tree2}, $F_t(T) \ge n_1$. In~\cite{k-Forcing}, it was shown that $F(T) \le n_1 - 1$ holds for all trees $T$. Consequently, $F_t(T) \ge F(T) + 1$, which establishes the desired inequality relating the total forcing number and forcing number of a non-trivial tree.
Suppose next that $F_t(T) = F(T) + 1$. Let $T'$ denote the trimmed tree of $T$; that is, $T' = \trim(T)$. By Lemma~\ref{lemTrees}, $F_t(T) = F_t(T')$, and by Lemma~\ref{forcingTreeLem}, $F(T) = F(T')$. Therefore, by supposition, $F_t(T') = F(T') + 1$.

Suppose that $T'$ has at least two strong support vertices, say $v$ and $w$. Let $S'$ be a minimum TF-set in $T'$, and so $|S'| = F_t(T')$. By Observation~\ref{ob:2}, the set $S'$ contains both $v$ and $w$, and all except possibly one leaf neighbor of each of $v$ and $w$ in $T'$. Let $v'$ and $w'$ be an arbitrary leaf neighbor of $v$ and $w$, respectively, that belongs to the set $S'$. We now consider the set $S = S' \setminus \{v,w\}$. We claim that $S$ is a forcing set of $T'$. In the first step of the forcing process starting with the initial set $S$, we play the vertex $v'$ which forces its (unique) neighbor $v$ to be colored. In the second step of the forcing process, we play the vertex $w'$ which forces its (unique) neighbor $w$ to be colored. At this stage of the forcing process, the resulting set of colored vertices is precisely the set $S'$, which is TF-set of $T'$ and therefore also a forcing set of $T'$. We now follow a sequence of played vertices in the total forcing process determined by the TF-set $S'$ of $T'$. In this way, all vertices of $V(T')$ are colored. Thus, the set $S$ is a forcing set of $T'$, implying that $F_t(T') - 1 = F(T') \le |S| = |S'| - 2 = F_t(T') - 2$, a contradiction. Therefore, $T'$ has at most one strong support vertex. If $T'$ has no strong support vertex, then $T$ is a path, and so $T \in \cH$. If $T'$ has exactly one strong support vertex, then $\trim(T) \cong K_{1,n_1}$ where $n_1 \ge 3$, and so $T \in \cH$. Hence, if $F_t(T) = F(T) + 1$, then $T \in \cH$.

It remains for us to prove that if $T \in \cH$, then $F_t(T) = F(T) + 1$. Let $T \in \cH$ have order~$n$ with $n_1$ leaves. If $T \cong P_n$, then it is well-known (and simple to observe) that $F(T) = 1 = n_1-1$. If $T \ncong P_n$, then $\trim(T) \cong K_{1,n-1}$ for some $n \ge 4$ and $F(T) = n - 2 = n_1-1$. In both cases, $F(T) = n_1-1$. Every tree in the family $\cH$ can be constructed from a path $P_2$ by applying a sequence of operations~$\cO_1$ and~$\cO_2$, and therefore belongs to the family $\cF$; that is, $\cH \subseteq \cF$. Hence since $T \in \cH$, we note that $T \in \cF$, implying by Lemma~\ref{l:Leaves2} that $F_t(T) = n_1$. Thus, $F_t(T) = F(T) + 1$.~\qed

\medskip
By Theorem~\ref{t:tree3}, all trees $T$ in the family~$\cH$ achieve equality in the inequality $F_t(T) \ge F(T) + 1$. We close by proving that the gap in this inequality can be made arbitrarily large.

\begin{prop}
\label{p:Gap}
For every integer $k \ge 1$, there exists a tree $T$ such that \[F_t(T) = F(T) + k.\]
\end{prop}
\proof Let $k \ge 1$ be an arbitrary integer, and let $T' \cong P_k$ be a path on $k$ vertices. Let $T$ be the graph obtained from $T'$ by adding two pendant edges to each vertex of $T'$. Thus, $T$ has order~$3k$ and every vertex in $V(T')$ is a strong support vertex of $T$. For each vertex of $V(T')$, select one of its leaf-neighbors and let $S$ denote the resulting set of $k$ leaves. The set $V(T) \setminus S$ is a TF-set of $T$, and so $F_t(T) \le |V(T)| - |S| = 2k$. Conversely, by Observation~\ref{ob:2}, $F_t(T) \ge 2k$. Consequently, $F_t(T) = 2k$. Moreover, the set $S$ is a forcing set of $T$, and so $F(T) \le k$. However, every forcing set of $T$ must contain at least one leaf neighbor of every vertex of $V(T')$ in $T$, implying that $F(T) \ge k$. Consequently, $F(T) = k$. Therefore, $F_t(T) - F(T) = k$.~\qed

\end{document}